\documentclass[12pt]{article}
\usepackage[american]{babel}
\usepackage{amsmath}
\usepackage{latexsym}
\usepackage{amssymb}
\usepackage{theorem}
\usepackage{diagrams}
\usepackage{mathrsfs}
\usepackage{calligra}
\DeclareMathAlphabet{\mathcalligra}{T1}{calligra}{m}{n}
\DeclareFontShape{T1}{calligra}{m}{n}{<->s*[1.2]callig15}{}
\theoremstyle{change}
\newtheorem{Thm}{Theorem}[section]

\newtheorem{Prop}[Thm]{Proposition}
\newtheorem{Lem}[Thm]{Lemma}
\newtheorem{Conj}[Thm]{Conjecture}
{\theorembodyfont{\rmfamily}
\newtheorem{Num}[Thm]{}

\newtheorem{Example}[Thm]{Example}
}

\renewcommand{\phi}{\varphi}

\newcommand{\proof}{\par\medskip\rm\emph{Proof. }}

\newcommand{\qed}{\ \hglue 0pt plus 1filll $\Box$}

\newcommand{\mapstoo}{\longmapsto}

\newcommand{\RR}{\mathbb{R}}

\newcommand{\SKIP}[1]{}

\renewcommand{\emptyset}{\varnothing}

\begin{document}

\title{{\bf Metric Properties of Euclidean Buildings}}
\author{Linus Kramer\thanks{The work by Essert, Hitzelberger, Kramer
and Weiss that we report on
was supported by the DFG through SPP 1154 \textit{Differentialgeometrie}, GK 
\textit{Analytische Topologie and Metageometrie}, SFB 478 and SFB 878.}}
\maketitle
\begin{abstract}
This is a survey on nondiscrete euclidean buildings, with a focus on
metric properties of these spaces.
\end{abstract}
Euclidean buildings are higher dimensional generalizations of trees.
Indeed, the euclidean product $X$ of two (leafless) metric trees $T_1$, $T_2$
is already a good ``toy example'' of a $2$-dimensional euclidean building.
The space $X$ contains lots of copies of the euclidean plane $\mathbb E^2$
and has at the same time a complicated local branching.

Euclidean building were invented by Jacques Tits in the seventies. Similarly
as in the case of spherical buildings, their definition was motivated by group theoretic
questions. While spherical buildings are by now a standard tool in the structure
theory of reductive algebraic groups over arbitrary fields,
euclidean buildings are important
for the advanced structure theory of reductive groups over fields with valuations.
In particular, they are very much linked to number theory and arithmetic geometry.

In the last 25 years, however, euclidean buildings have also become
important in geometry.
This is due to the fact that euclidean buildings are spaces of nonpositive
curvature. But more is true. Together with the Riemannian symmetric spaces of
nonpositive curvature, euclidean buildings could be called the islands of high symmetry
in the world of CAT$(0)$ spaces. This claim will be made more precise below.
Almost inevitably, questions about symmetry, rigidity, or higher rank
for CAT$(0)$ spaces lead to these geometries.

I thank Alexander Lytchak, Koen Struyve and Richard Weiss for helpful comments.

\section{The definition of euclidean buildings}

We first recall Tits' definition of a euclidean building. For more details,
proofs and further results see
Tits \cite{TitsComo}, Kleiner and Leeb \cite{KL}, Kramer and Weiss \cite{KrWe} and in particular 
Parreau \cite{Par}.
(The axioms used by Kleiner and Leeb \cite{KL}
look somewhat different from Tits' definition. They were shown to be
equivalent to Tits' by Parreau.)
\begin{Num}\textbf{Euclidean buildings}\label{EBAxioms}\\
Let $W$ be a spherical Coxeter group acting in its natural
orthogonal representation on euclidean space $\mathbb E^m$. We call the
semidirect product $W\RR^m$ of $W$ and $(\RR^m,+)$ the 
\textit{affine Weyl group}.
From the reflection hyperplanes of $W$ we
obtain a decomposition of $\mathbb R^m$ into \emph{walls}, \emph{half spaces},
\emph{Weyl chambers} (a Weyl chamber
is a fundamental domain for $W$---these are Tits' \emph{chambres vectorielles}) and 
\emph{Weyl simplices} (Tits' \emph{facettes vectorielles}).
The $W\RR^m$-translates of these in $\mathbb E^m$ are also called walls,
half spaces and Weyl chambers.

Let now $X$ be a metric space. A \emph{chart} is an isometric
embedding $\phi:\mathbb E^m\rTo X$, and its image is called an (affine) \emph{apartment}.
We call two charts $\phi,\psi$
\emph{$W$-compatible} if $Y=\phi^{-1}(\psi(\mathbb E^m))$ is convex 
(in the Euclidean sense) and
if there is an element $w\in W\RR^m$ such that
$\psi\circ w|_Y=\phi|_Y$ (this condition is void if $Y=\emptyset$). 
We call a metric space $X$ together with a collection $\mathcal A$ of
charts a \emph{Euclidean building} if it has the following
five properties.
\begin{enumerate}
\item[(A1)] For all $\phi\in\mathcal A$ and $w\in W\RR^m$, the composition
$\phi\circ w$ is in $\mathcal A$.
\item[(A2)] The charts are $W$-compatible.
\item[(A3)] Any two points $p,q\in X$ are contained in some affine apartment.
\end{enumerate}
The charts allow us to map Weyl chambers, walls and half spaces into $X$;
their images are also called Weyl chambers, walls and half spaces. The first
three axioms guarantee that these notions are coordinate independent.
\begin{enumerate}
\item[(A4)] If $C,D\subseteq X$ are Weyl chambers, then there is
an affine apartment $A$ such that the intersections
$A\cap C$ and $A\cap D$ contain Weyl chambers.
\item[(A5')] For every apartment $A\subseteq X$ and every $p\in A$ there is
a $1$-Lipschitz retraction $h:X\rTo A$ with $h^{-1}(p)=\{p\}$.
\end{enumerate}
Condition (A5') may be replaced by the following condition:
\begin{enumerate}
\item[(A5)] If $A_1,A_2,A_3$ are affine apartments which intersect pairwise
in half spaces, then $A_1\cap A_2\cap A_3\neq\emptyset$.
\end{enumerate}
See \cite{Par} for a thorough discussion of different sets of axioms.

Let $p$ be a point in the apartment $A$.
Axiom (A5') yields a $1$-Lipschitz map $v_p:X\rTo\mathbb E^m/W$ as follows.
We identify $(A,p)$ by means of a coordinate chart $\phi$ with $(\mathbb E^m,0)$,
and then we quotient out the $W$-action. The resulting vector $v_p(q)$ is called the
\textit{vector distance} between $p$ and~$q$.
\end{Num}
The definition of a euclidean building that we use here is the metric,
``nondiscrete'' version. It appeared implicitly in \cite{BruTi}, and
in detail in \cite{TitsComo}. There is also the (older) notion of a simplicial
affine building; see \cite{AB} and in particular \cite{WeAff}.
The geometric realization of such
a combinatorial affine building is always a euclidean building in our
sense (but not vice versa); see \S11.2 in \cite{AB}.

An important invariant of a euclidean building is its \textit{spherical building
at infinity}, $\partial_{\mathcal A}X$. This is a combinatorial simplicial
complex which is defined as follows. The simplices are equivalence
classes of Weyl simplices. Two Weyl simplices are considered to be equivalent
if they have finite Hausdorff distance. It turns out that $\partial_{\mathcal A}X$
is a (weak) \textit{spherical building} of dimension $m-1$ (resp., of rank $m$)
\cite[\S1.5]{Par}.
We refer to \cite{TitsLNM} and \cite{AB}
for the definition of a simplicial spherical building.

Another important fact is that a euclidean building is always a CAT$(0)$ space;
see \cite[\S2.3]{Par}.
This was first observed by Bruhat and Tits in \cite{BruTi}, where they
proved the CN-inequality.

\textit{Automorphisms} of euclidean buildings 
are defined in the obvious way; they are bijections
which preserve the charts in the given atlas. Clearly, every automorphism is
an isometry of $X$.

\section{Basic and less basic properties}
We assume that $X$ is a euclidean building with $m$-dimensional apartments.
First of all, we remark that the atlas $\mathcal A$ is by 
no means unique. However, Parreau proved that there is a unique
maximal atlas $\mathcal A_{max}$ containing $\mathcal A$; see \cite[\S2.6]{Par}.
The apartments in the maximal atlas have a simple characterization.
\begin{Thm}[{{\cite[2.6]{Par} and \cite[\S4.6]{KL}}}]
Let $X$ be a euclidean building with $m$-dimensional apartments. Suppose that
$F\subseteq X$ is a subspace isometric to some $\mathbb E^\ell$. Then there exists
an apartment in the maximal atlas containing $F$. In particular, $\ell\leq m$
and the apartments of the maximal atlas are precisely the maximal flats in $X$.
\end{Thm}
The metric realization of  the spherical building
$\partial_{\mathcal A_{max}}X$ can be identified in a canonical
way with the Tits boundary of the CAT$(0)$ space $X$.
We remark that the dimension $m$ of the apartments coincides with the covering
dimension of $X$ as a topological space; see  \cite[Prop.~3.3]{LaS} or \cite[Thm.~B]{KramerLocal}.
Moreover, $X$ is an AR (an absolute retract for the class of metric spaces).

Next, we note that the Weyl group may be ``too big'': there might be types of walls
which never appear as branchings between apartments. A wall $M$ in a euclidean
building $X$ is called \textit{thick} if it can be written as the intersection
of three apartments. We call a point $p\in X$ \textit{thick} if every wall
passing through $p$ is thick. Now we can make the statement about the Weyl
group being too big more precise: if $X$ contains no thick points, then
there is a (unique) euclidean building $X_{th}$ (with a smaller Weyl group) containing
a thick point, and $X$ is a euclidean product
$X\cong X_{th}\times\mathbb E^k$, for some $k\geq 0$; 
see \cite[\S4.9]{KL} and \cite[\S10]{KrWe}.
For the thick part $X_{th}$, there is the following trichotomy.
\begin{Prop}[{{\cite[\S10]{KrWe}}}]
\label{Trichotomy}
Let $X$ be an irreducible euclidean building of dimension $m\geq 2$
containing a thick point. Then there are the following three possibilities.

\medskip\noindent
(I) There is a unique thick point which is contained in every affine apartment
of $X$. In this case $X$ is a euclidean cone over a spherical building.\\
(II) The set of thick points is a closed, discrete and cobounded
subset in $X$ and in  every apartment of $X$. Then $X$ is the geometric
realization of a simplicial affine building.\\
(III) The set of thick points is dense in $X$ and in every apartment of $X$.
\end{Prop}
A simplicial affine building (type II) is called \textit{thick} if every vertex of the
simplicial structure is thick.

There are many $2$-dimensional euclidean buildings. In fact, there are ``free
constructions'' which show that it is impossible to classify these spaces.
In higher dimensions, the picture is completely different.
We call a euclidean building $X$ a \textit{Bruhat-Tits building}
if the spherical building at infinity is a \textit{Moufang building};
see \cite{WeSph}. Roughly speaking, the Moufang condition says that
there are certain automorphisms, called \textit{root automorphisms},
that fix a large subset pointwise, and yet act transitively on another
subset.
The following deep result is again due to Tits \cite{TitsComo}.
\begin{Thm}
\label{BTB}
Let $X$ be an irreducible euclidean building of dimension $m\geq 3$
containing a thick point. Then $\partial_{\mathcal A}X$ is Moufang,
and all root automorphisms of $\partial_{\mathcal A}X$ extend to
isometries of $X$. In particular, the isometry group of $X$ is
transitive on the apartments of $X$.
\end{Thm}
Tits' article \cite{TitsComo} contains in fact a complete classification
of these buildings in terms of algebraic data. 
We remark that if a Bruhat-Tits building is not of type (I),
the group generated by the root automorphisms acts with cobounded orbits on $X$.

It is by no means clear that every combinatorial automorphism of
$\partial_{\mathcal A}X$ extends to an isometry of $X$. Surprisingly, the
following is true.
\begin{Thm}[{{\cite[27.6]{WeAff}}}]
\label{CompleteFields}
Let $X$ be a thick irreducible simplicial Bruhat-Tits building of dimension $m\geq 2$.
Then every automorphism of $\partial_{\mathcal A_{max}}X$ extends to an isometry
of $X$. Moreover, $\partial_{\mathcal A_{max}}X$ determines $X$ up to isomorphism.
\end{Thm}
The proof depends on the purely algebraic fact that a field admits at most
one discrete complete valuation. It would be interesting to have a geometric
proof for this. More generally, there is the following open problem.
\begin{Num}\textbf{Question}
\em Is a thick irreducible simplicial affine building of dimension $m\geq 2$
uniquely determined by the spherical building $\partial_{\mathcal A_{max}}X$?
Does every combinatorial automorphism of $\partial_{\mathcal A_{max}}X$ extend
to $X$?
\end{Num}
The answer is negative if $X$ is not assumed to be simplicial.
For locally finite simplicial thick irreducible affine buildings, the answer
is positive \cite{LeebHabil}.
%
%
%
%

\section{Characterizations}

The following very general characterization of locally finite (simplicial) euclidean buildings
is due to Kleiner.
\begin{Thm}
Let $X$ be a locally compact CAT$(0)$ space of dimension $m$. Suppose that any two
points $x,y\in X$ are contained in some flat $A\cong\mathbb E^m$. Then $X$ is
a euclidean building.
\end{Thm}
This result was not published by Kleiner; a proof was given by Balser and Lytchak
in \cite[Cor.~1.7]{BL1}.
The dimension may be taken to be the covering dimension; since $X$ is locally
compact, the covering dimension coincides with Kleiner's geometric dimension \cite{Kl}.
The following example shows that local compactness is crucial.
\begin{Example}
Let $\Gamma_n$, for $n\geq 3$, be a family of thick generalized $n$-gons
($1$-dimensional spherical buildings whose Weyl group is dihedral of order $2n$).
Such generalized $n$-gons exist by Tits' free construction \cite{TitsFree},
see also \cite{TentHom}.
Let $X_n$ be the euclidean cone over $\Gamma_n$, with cone point $o_n$. Then $X_n$ is
a $2$-dimensional euclidean building with precisely one thick point. Now consider the
asymptotic cone (or ultralimit)  $X$ over the family $\{(X_n,o_n)\mid n\geq 3\}$ (with respect to a constant
scaling sequence and a nonprincipal ultrafilter $\mu$ on the index set
$\mathbb N_{\geq 3}$).
Then $X$ is a complete CAT$(0)$ space. Any two points in $X$
are contained in some copy of $\mathbb E^2$. 
The ``spherical Weyl group'' $W$ that describes the
transition functions between these ``apartments'' is, however, the orthogonal group
$W=\mathrm O(2)$. Using a similar argument as in 
\cite[\S7]{KramerLocal} (or by Kleiner's results in \cite{Kl}, 
see also Lytchak \cite[11.3]{LytRig}) one can show that $X$ is $2$-dimensional. 
But $X$ is certainly not a euclidean building.
Instead of the cones $X_n$, one could also use the euclidean buildings constructed
recently by Berenstein and Kapovich \cite{BeKa}
in order to get a more interesting asymptotic cone $X$.
\end{Example}
In a somewhat more combinatorial setting, there is the following result of
Charney and Lytchak. A CAT$(0)$ space $X$ has the \textit{discrete extension property}
if for every geodesic $\gamma=[a,b]\subseteq X$, the set of the directions of
geodesics extending $\gamma$ beyond $b$ is nonempty and discrete.
\begin{Thm}
Let $X$ be a CAT$(0)$ space of dimension $m\geq 2$
which is a piecewise euclidean cell complex.
If $X$ has the discrete extension property, then $X$ is a euclidean building.
\end{Thm}
We remark that a locally compact euclidean building always admits a
euclidean cell structure. This is not true for general euclidean buildings.
Finally, we should mention here the following result by Leeb \cite{LeebHabil}.
\begin{Thm}
Let $X$ be a locally compact CAT$(0)$ space with extendible geodesics. If the
Tits boundary of $X$ is an irreducible spherical building of rank at least $2$,
then $X$ is either a Riemannian symmetric space of noncompact type or
a simplicial euclidean building.
\end{Thm}

\section{Isometries and automorphisms}

If $X$ is a euclidean building containing a thick point, then an isometry of
$X$ is almost the same as an automorphism.
\begin{Thm}[{{\cite[\S4]{Par}}}]
Let $g$ be an isometry of a euclidean building $X$. Assume that $X$ contains a
thick point. Then there exists an element $\gamma\in \mathrm{Nor}_{\mathrm O(m)}(W)$
such that $g\circ\psi\circ \gamma\in\mathcal{A}_{max}$ holds for all
$\psi\in\mathcal{A}_{max}$.
\end{Thm}
Such a map $\gamma$ induces a diagram automorphism of the Coxeter group $W$;
one also calls such a $g$ a \textit{non-type-preserving automorphism}. 

Suppose that $g$ is an isometry of a metric space $(X,d)$. The \textit{displacement
function} of $g$ is the nonnegative real function $d_g:x\mapstoo d(x,g(x))$. 
The infimum of $d_g(X)$ is the \textit{translation length} $l_g$ of $g$.
We call an isometry $g$

\medskip
\textit{elliptic} if $g$ has a fixed point.

\textit{hyperbolic} if $d_g$ attains a positive minimum.

\textit{parabolic} if $d_g$  does not attain its minimum.

\medskip\noindent
If $X$ is a Riemannian symmetric space of nonpositive curvature, then all three types of
isometries appear in the isometry group. This is not true for euclidean buildings.
\begin{Thm}[{{\cite[\S4]{Par}}}]
Let $g$ be an isometry of a complete euclidean building $X$ containing a thick
point. Then $g$ is either elliptic or hyperbolic.
\end{Thm}
(Struyve informed me that he can prove this also for noncomplete
euclidean buildings.)
The next result was proved by Rapoport and Zink \cite{RZ} for the Bruhat-Tits building
of $\mathrm{GL}_n$ over a field with discrete valuation, and then extended using
Landvogt's Embedding Theorem to other Bruhat-Tits buildings. However, there is
a much simpler proof using CAT$(0)$ geometry, which applies to all euclidean
buildings, cp.~\cite{Rou}---the author found a somewhat simpler proof (unpublished).
We put $X_r=\{q\in X\mid d_g(q)\leq r\}$. These sublevel sets form a
filtration of $X$ by convex sets.
\begin{Thm}
Let $g$ be an isometry of a complete euclidean building $X$ containing a thick point.
There exists a positive constant $c>0$ (depending only on the Weyl group $W$) such that
the following holds.
If $p$ is a point with $d(p,X_r)=t>0$, then 
\[
c\cdot t+r\leq d_g(p)\leq2t+r.
\]
\end{Thm}
The second inequality is trivial, the interesting fact is the lower estimate.
We finally note the following (completely elementary) fact.
\begin{Lem}
Let $g$ be a nontrivial isometry of a euclidean building $X$ containing a thick point.
Then $\sup d_g(X)=\infty$.

\proof
Suppose $r=\sup d_g(X)<\infty$. If $A$ is an apartment in $X$, then
$g(A)$ has Hausdorff distance at most $r$ from $A$. Then $A$ and $g(A)$
have the same boundary at infinity. By \cite[p.~10]{Par}, $A=g(A)$. Thus
$g$ fixes all apartments setwise, and therefore all thick walls and thick points.
Since every apartment contains a thick point, $g$ fixes every apartment pointwise.
Thus $g=id_X$.
\qed 
\end{Lem}
We end this section with some remarks on noncomplete euclidean buildings.
Struyve recently proved the following generalization of the Bruhat-Tits
Fixed Point Theorem. If a finitely generated group acts acts isometrically and
with bounded orbits on a euclidean building, then it has a fixed point \cite{Struyve}.
Moreover, he showed that the main rigidity results in \cite{KrWe} also hold if the
completeness assumptions on the euclidean buildings are dropped (unpublished).
Finally, he and Martin, Schillewaert and Steinke extended results in \cite{BruTi}
about noncomplete Bruhat-Tits buildings, by giving algebraic conditions
on the underlying fields (unpublished).

\section{Kostant convexity}

We first recall the statement of Kostant's Convexity Theorem \cite{Kostant}
for Riemannian symmetric spaces.
Let $G$ be a simple
noncompact Lie group with Iwasawa decomposition $G=K\!AU$. The group $K$ is maximal
compact, $A$ is diagonizable, and $U$ is unipotent. The group 
$W=\mathrm{Nor}_K(A)/\mathrm{Cen}_K(A)$ is the associated Weyl group.

The solvable group $AU$ acts
regularly on the Riemannian symmetric space $X=G/K$. Let $o\in X$ denote the
point stabilized by $K$. The $A$-orbit $E=A(o)\subseteq X$
is a maximal flat in $X$.
The projection $AU\rTo AU/U\cong A$ induces a natural $1$-Lipschitz
map $\rho_U:X\rTo E$ which we call the
\textit{Iwasawa projection}. Let $p\in E$. The Convexity Theorem says that 
\[
\rho_U(K(p))=\mathrm{conv}(W(p)),
\]
the image of the $K$-orbit of $p$ in $X$ under the
Iwasawa projection is the convex hull of the $W$-orbit of $p$ in $E$.

Geometrically, the Iwasawa decomposition can also be described as follows.
The group $U$ determines a chamber $C$ of the spherical building at infinity
of $X$. The maximal flats in $X$ containing $C$ in their boundary form a foliation
of $X$. The Iwasawa projection  identifies each leaf by means of the $U$-action with
the leaf $A(o)$.

Suppose now that $X$ is a euclidean building and that $C$ is a chamber
at infinity. We fix an apartment $E\subseteq X$ containing $C$ in its
boundary. If $E'\subseteq X$ is any other apartment containing
$C$ in its boundary, then $E\cap E'$ contains a Weyl chamber representing $C$.
Thus, there is a canonical isometry $E'\rTo E$ fixing $E\cap E'$ pointwise.
These isometries fit together to a $1$-Lipschitz retraction 
$\rho_C:X\rTo E$. Suppose now that $o,p\in E$ are special vertices. 
(A vertex $p\in E$ is called special if the reflections along the thick walls 
in $E$ passing through $p$ generate the spherical Weyl group $W$.)
Let
$S\subseteq X$ denote the set of all special vertices in $X$ that have the
same vector distance from $o$ as $p$. 
This set $S$ corresponds to the
orbit $K(p)$ in the Riemannian symmetric case. If the euclidean building $X$
happens to be a Bruhat-Tits building, then $S$ is indeed the $K$-orbit of
$p$, where $K$ is the stabilizer of $o$.
The following result was proved by Hitzelberger \cite{HitzCvx} in 2007.
\begin{Thm}
Suppose that $X$ is a thick simplicial euclidean building. With the same notation
as above, suppose that $o,p$ are special vertices (see \cite{BruTi} for the
definition of a special vertex). Then 
\[
\rho_C(S)=\{q\in\mathrm{conv}(W(p))\mid q\text{ has the same type as }p\},
\]
the image of $S$ is the set of all vertices in $E$ which are in the
convex hull of the $W$-orbit of $p$ in $E$ and have the same type as $p$.
\end{Thm}
This result had been announced by Silberger \cite{Sil} for the special case that
$X$ is the Bruhat-Tits building of a simple $p$-adic algebraic group
(but the proof, which relied on a case-by-case analysis, was never
published). The difficult part of the proof is to show that the map
is onto. For the special case of
Bruhat-Tits buildings, the theorem may be restated as
a fact about intersections of certain double cosets in the group.
The result was recently extended by Hitzelberger to
general euclidean buildings \cite{HitzCvx2}.

\section{Rigidity}

We first recall some notions from coarse geometry \cite{Roe}.
A map $f:X\rTo Y$ between metric spaces is called
\emph{controlled} if there is a monotonic real function
$\rho:\RR_{\geq 0}\rTo\RR_{\geq 0}$ such that
\[
d_Y(f(x),f(y))\leq\rho(d_X(x,y))
\]
holds for all $x,y\in X$. If in addition the preimage of every
bounded set is bounded, then $f$ is called a \emph{coarse map}.
Neither $f$ nor $\rho$ is required to be continuous.
Note that the image of a bounded set under a controlled map is
bounded.
Two maps $g,f:X\pile{\rTo \\ \rTo} Y$ between metric spaces have \emph{finite distance}
if the set $\{d_Y(f(x),g(x))\mid x\in X\}$ is bounded. This is an
equivalence relation which
leads to the \emph{coarse metric category} whose
objects are metric
spaces and whose morphisms are equivalence classes of coarse maps.
A \emph{coarse equivalence} is an isomorphism in this category.
We remark that a coarse equivalence between geodesic metric spaces is
the same as a quasi-isometric equivalence.

Prasad proved in 1978
the following analog of Mostow's Rigidity Theorem.
\begin{Thm}[\cite{Pra}]
Let $X$ and $Y$ be thick simplicial, irreducible and locally finite Bruhat-Tits
buildings of rank at least $2$. Suppose that a group $\Gamma$ acts cocompactly and
properly discontinuously on both spaces. Then there is a $\Gamma$-equivariant
simplicial isomorphism between $X$ and $Y$.
\end{Thm}
The group $\Gamma$ appearing in Prasad's Theorem is finitely presentable.
From the $\Gamma$-action, one obtains a $\Gamma$-equivariant coarse equivalence
$f:X\rTo Y$ which plays a crucial role in the proof.
About twenty years later, Kleiner and Leeb \cite{KL} proved the following generalization
of Prasad's Theorem.
\begin{Thm}[\cite{KL}]
Let $X$ and $Y$ be complete Bruhat-Tits buildings whose de Rham factors
all have rank at least $2$. Suppose that $f:X\rTo Y$ is a coarse
equivalence. Then there is an isometry $\bar f:X\rTo Y$ 
(possibly after rescaling the metrics on the de Rham factors of $Y$)
which has finite distance from $f$.
\end{Thm}
The strategy of their proof is roughly as follows. Using \textit{asymptotic cones},
Kleiner and Leeb show that the $f$-image of a maximal flat $E\subseteq X$
has finite Hausdorff distance from a (necessarily unique) maximal flat
$E'\subseteq Y$.
This fact is then used to set up a one-to-one correspondence between
the maximal bounded subgroups of the isometry groups of the two
Bruhat-Tits buildings. The maximal bounded subgroups, in turn, correspond
to (certain) points in the buildings. In this way, they construct an
equivariant isometry.

Weiss and the author proved in 2009 a more general result which
is valid for all euclidean buildings.
\begin{Thm}[{{\cite[Thm.~III]{KrWe}}}]
\label{KW1}
Let $X$ and $Y$ be complete euclidean buildings containing thick points, 
and without rank $1$ de Rham factors.
Suppose that $f:X\rTo Y$ is a coarse equivalence.
Then there is an isometry $\bar f:X\rTo Y$. If no de Rham factor
of $X$ is a euclidean cone, then $f$ has finite distance from $\bar f$.
\end{Thm}
The proof relies, among other things, on the following result about trees.
\begin{Thm}[{{\cite[Thm.~I]{KrWe}}}]
\label{KW2}
Let $T,T'$ be two complete $\mathbb R$-trees without leaves. Suppose that a group
$G$ acts isometrically
on both trees, and that this action is $2$-transitive on the
ends. Suppose that $f:T\rTo T'$ is a coarse equivalence whose
induced boundary map $\partial T\rTo\partial T'$ is $G$-equivariant.
Then $T$ and $T'$ are $G$-equivariantly isometric.
\end{Thm}
The proof of \ref{KW1} proceeds roughly as follows. The first step
is a result due to Kleiner and Leeb which was already mentioned: 
the $f$-image of an apartment $E\subseteq X$
has finite Hausdorff distance from a
(unique) apartment $E'\subseteq Y$.
But then we follow a different line. We show
directly that $f$ induces a combinatorial isomorphism $f_*$
between the Tits boundaries of $X$ and $Y$.
(For the case of
simplicial Bruhat-Tits buildings, this implies
by \ref{CompleteFields} already that $X$ and $Y$ are combinatorially isomorphic.)
Next, we show that we obtain a coarse bijection between the
so-called \textit{wall trees} of $X$ and $Y$. Since these trees have large
holonomy groups, we may apply \ref{KW2}. In this way we get equivariant
isomorphisms between
the wall trees, and thus, by Tits \cite{TitsComo}, an isometry
between the euclidean buildings. We remark that the main results in \cite{TitsComo}
also enter as important ingredients into the proof of \cite{KL}.

\section{Locally compact Bruhat-Tits buildings}

In the mid-nineties, Grundh\"ofer, Knarr and the author completed
the classification of all compact connected spherical buildings
admitting a chamber transitive automorphism group. Such buildings
arise for example as boundaries of Riemannian symmetric spaces.
The proof and the method of the classification
built on earlier work by Salzmann, L\"owen,
Burns and Spatzier. Briefly, it may be stated as follows.
\begin{Thm}[\cite{GKK1,GKK2,kramerhabil}]
\label{GKK}
Let $B$ be a compact spherical building (in the sense of \cite{burns})
without rank $1$ factors.
Suppose that $B$ is (locally) connected and admits a chamber
transitive group of continuous automorphisms. Then $B$ is the
Tits boundary of a Riemannian symmetric space of noncompact type.
\end{Thm}
There should be an analog of this result, corresponding to
the boundaries of locally compact euclidean buildings.
The following conjecture is wide open (even for buildings of type
$A_2$, i.e. compact projective planes).
\begin{Conj}
Let $B$ be a compact spherical building (in the sense of \cite{burns})
without rank $1$ factors.
Suppose that $B$ is totally disconnected and admits a chamber
transitive groups of continuous automorphisms. Then $B$ is the
Tits boundary of a locally finite simplicial Bruhat-Tits
building.
\end{Conj}
The problem is that in comparison to \ref{GKK},
no homotopy theory is available. Presently, a proof of this
conjecture seems to be out of reach.
Assuming the Moufang property, we showed however the following.
\begin{Thm}[\cite{GKVW}]
\label{GKVW}
Let $B$ be a compact spherical building (in the sense of \cite{burns})
without rank $1$ factors.
Suppose that $B$ is totally disconnected and Moufang.
Then $B$ is the
Tits boundary of a locally finite simplicial Bruhat-Tits
building.
\end{Thm}
We recall that the Moufang property is automatically satisfied
if all irreducible factors of $B$ have rank at least $3$;
see \cite{TitsLNM} and \cite{WeSph}.
The proof of \ref{GKVW} relies very much on the classification
of spherical Moufang buildings due to Tits and Weiss.

\section{Lattices}

Let $X$ be a complete and locally compact CAT$(0)$ space and let
$\Gamma$ be a group of isometries. We call $\Gamma$ a 
\textit{uniform lattice}
if $\Gamma$ acts properly discontinuously and cocompactly on $X$
(such groups are also called CAT$(0)$ groups).
Borel's Density Theorem says that Riemannian symmetric spaces
of noncompact type admit (many) uniform lattices. 
Such a uniform lattice is always finitely presentable. However, very
few presentations of lattices are known.
Essert observed the following correspondence between uniform lattices
acting regularly on the $1$-simplices of a given type of a $2$-dimensional
locally finite simplicial euclidean building and Singer groups.
A \textit{Singer group} is a subgroup of the automorphism group a finite
generalized polygon (a $1$-dimensional spherical building) which acts
regularly on the vertices of a given type. Singer groups are studied by
finite geometers and group theorists, and quite a few constructions are
known. Essert showed that from
a collection of Singer groups, one can construct a $2$-dimensional
\textit{complex of groups} which unfolds to a lattice $\Gamma$
acting on such a $2$-dimensional euclidean building.
Specific examples are presentations such as
\[
\langle a,b,c\mid a^7=b^7=c^7=abc=a^3b^3c^3=1\rangle
\]
or
\[
\langle a,b,c\mid a^{13}=b^{13}=c^{13}=ab^3c^9=a^3b^9c=a^9bc^3=1\rangle.
\]
These explicit representations allow,
for example, to compute the group homology of the lattices. It is clear that
``most'' of the buildings $X$ that he constructed in this way are
``exotic'', i.e. they are not Bruhat-Tits buildings. There are presently
many open questions about these lattices $\Gamma$, eg. about commensurabilty,
quasi-isometric type, or the covolume. We refer to \cite{EssertLattices} for
details and more results.

\section{Noncrystallographic Weyl groups}

The Weyl group of a Bruhat-Tits building arising from a reductive
algebraic group over a field with valuation is always crystallographic.
Also, the Weyl group of a simplicial euclidean building is
necessarily a crystallographic group. But in the definition of a euclidean
buildings, there is no reason to assume that $W$ satisfies the
crystallographic condition. It was remarked (without giving details)
by Tits \cite{TitsComo} that there are Bruhat-Tits buildings with
non-crystallographic Weyl groups. An explicit construction of such
euclidean buildings, defined over certain, very special fields, can be found in
\cite{HKW}. Their Weyl groups are dihedral groups of order $16$,
and their Tits boundaries are so-called Moufang generalized octagons.

In a completely different way, Berenstein and Kapovich
constructed ``wild'' $2$-dimensional euclidean buildings
whose Weyl groups are dihedral groups of arbitrary order
\cite{BeKa}. It would be interesting to see if the
construction can be done in such a way that it yields
highly transitive automorphism groups, 
as was the case for the $1$-dimensional spherical buildings
constructed by Tent in\cite{TentHom}.

\raggedright
Linus Kramer\\
Mathematisches Institut, 
Universit\"at M\"unster,
Einsteinstr. 62,
48149 M\"unster,
Germany\\
\makeatletter
e-mail: {\tt linus.kramer{@}uni-muenster.de}

\begin{thebibliography}{88}

 
\bibitem{AB}
P. Abramenko\ and\ K. S. Brown, {\it Buildings}, Graduate Texts in Mathematics, 248, Springer, New York, 2008. 

 
\bibitem{BGS}
W. Ballmann, M. Gromov\ and\ V. Schroeder, {\it Manifolds of nonpositive curvature}, Progress in Mathematics, 61, Birkh\"auser Boston, Boston, MA, 1985. 

 
\bibitem{BL1}
A. Balser\ and\ A. Lytchak, Building-like spaces, J. Math. Kyoto Univ. {\bf 46} (2006), no.~4, 789--804. 

\bibitem{BeKa}A. Berenstein and M. Kapovich,
Affine buildings for dihedral groups. 
To appear in: Geom. Dedicata.

 
\bibitem{BH}
M. R. Bridson\ and\ A. Haefliger, {\it Metric spaces of non-positive curvature}, Grundlehren der Mathematischen Wissenschaften, 319, Springer, Berlin, 1999. 

 
\bibitem{BruTi}
F. Bruhat\ and\ J. Tits, Groupes r\'eductifs sur un corps local, Inst. Hautes \'Etudes Sci. Publ. Math. No. 41 (1972), 5--251. 

%

 
\bibitem{burns}
K. Burns\ and\ R. Spatzier, On topological Tits buildings and their classification, Inst. Hautes \'Etudes Sci. Publ. Math. No. 65 (1987), 5--34. 

 
\bibitem{CL}
R. Charney\ and\ A. Lytchak, Metric characterizations of spherical and Euclidean buildings, Geom. Topol. {\bf 5} (2001), 521--550. 


\bibitem{EssertLattices}
J. Essert,
A geometric construction of panel-regular lattices in buildings of types $\tilde A_2$ 
and $\tilde C_2$.
Preprint 2009, submitted.
aXiv:0908.2713

 
\bibitem{GKK1}
T. Grundh\"ofer, N. Knarr\ and\ L. Kramer, Flag-homogeneous compact connected polygons, Geom. Dedicata {\bf 55} (1995), no.~1, 95--114. 

 
\bibitem{GKK2}
T. Grundh\"ofer, N. Knarr\ and\ L. Kramer, Flag-homogeneous compact connected polygons. II, Geom. Dedicata {\bf 83} (2000), no.~1-3, 1--29. 

\bibitem{GKVW}
T. Grundh\"ofer, L. Kramer, H. Van Maldeghem and R. M. Weiss,
Compact totally disconnected Moufang buildings.
To appear in: Tohoku Math. Journal.

\bibitem{HitzCvx}
P. Hitzelberger, Kostant convexity for affine buildings. Forum Math. {\bf22} (2010), 959--971.

\bibitem{HitzCvx2}
P. Hitzelberger, Non-discrete affine buildings and convexity.
Advances in Mathematics {\bf 227} (2011), no.~1, 210--244.

 
\bibitem{HKW}
P. Hitzelberger, L. Kramer\ and\ R. M. Weiss, Non-discrete Euclidean buildings for the Ree and Suzuki groups, Amer. J. Math. {\bf 132} (2010), no.~4, 1113--1152. 

 
\bibitem{Kl}
B. Kleiner, The local structure of length spaces with curvature bounded above, Math. Z. {\bf 231} (1999), no.~3, 409--456. 

 
\bibitem{KL}
B. Kleiner\ and\ B. Leeb, Rigidity of quasi-isometries for symmetric spaces and Euclidean buildings, Inst. Hautes \'Etudes Sci. Publ. Math. No. 86 (1997), 115--197 (1998). 

 
\bibitem{Kostant}
B. Kostant, On convexity, the Weyl group and the Iwasawa decomposition, Ann. Sci. \'Ecole Norm. Sup. (4) {\bf 6} (1973), 413--455 (1974). 

 
\bibitem{kramerhabil}
L. Kramer, Homogeneous spaces, Tits buildings, and isoparametric hypersurfaces, Mem. Amer. Math. Soc. {\bf 158} (2002), no.~752, xvi+114 pp. 

\bibitem{KramerLocal}
L. Kramer,
On the local structure and the homology of CAT$(\kappa)$ spaces and euclidean buildings.
Adv. Geom. 11 (2011), 347–369.
 
\bibitem{KrTe}
L. Kramer\ and\ K. Tent, Asymptotic cones and ultrapowers of Lie groups, Bull. Symbolic Logic {\bf 10} (2004), no.~2, 175--185. 

\bibitem{KrWe}
L. Kramer\ and\ R. M. Weiss, 
Coarse rigidity of euclidean buildings. Preprint 2009, submitted. arXiv:0902.1332

\bibitem{LaS}
U. Lang\ and\ T. Schlichenmaier, Nagata dimension, quasisymmetric embeddings, and Lip\-schitz extensions, Int. Math. Res. Not. {\bf 2005}, no.~58, 3625--3655. 

\bibitem{LeebHabil}
B. Leeb, {\it A characterization of irreducible symmetric spaces and Euclidean buildings of higher rank by their asymptotic geometry}, 
Bonner Mathematische Schriften, 326. Univ. Bonn, Bonn, 2000. 

\bibitem{LytRig}
A. Lytchak, Rigidity of spherical buildings and joins, Geom. Funct. Anal. {\bf 15} (2005), no.~3, 720--752. 

\bibitem{Par}
A. Parreau, Immeubles affines: construction par les normes et \'etude des isom\'etries, in {\it Crystallographic groups and their generalizations (Kortrijk, 1999)}, 263--302, Contemp. Math., 262 Amer. Math. Soc., Providence, RI. 


 
\bibitem{Pra}
G. Prasad, Lattices in semisimple groups over local fields, in {\it Studies in algebra and number theory}, 285--356, Adv. in Math. Suppl. Stud., 6 Academic Press, New York. 

\bibitem{RZ} 
M. Rapoport and T. Zink,
A finiteness theorem in the Bruhat-Tits building:
an application of Landvogt's embedding theorem.
Indag. Math. (N.S.)  10  (1999),  no. 3, 449--458.
 
\bibitem{Roe}
J. Roe, {\it Lectures on coarse geometry}, University Lecture Series, 31, Amer. Math. Soc., Providence, RI, 2003. 

\bibitem{Rou} G. Rousseau,
Exercices m\'etriques immobiliers. 
Indag. Math. (N.S.)  12  (2001),  no. 3, 383--405.


\bibitem{Sil}
A. J. Silberger, Convexity for a simply connected $p$-adic group, Bull. Amer. Math. Soc. {\bf 81} (1975), no.~5, 910--912. 

 

\bibitem{Struyve}
K. Struyve,
(Non-)completeness of $\mathbb R$-buildings and fixed point theorems.
Groups, Geometry and Dynamics {\bf 5} (2011), no.~1, 177--188.
 
\bibitem{TentHom}
K. Tent, Very homogeneous generalized $n$-gons of finite Morley rank, J. London Math. Soc. (2) {\bf 62} (2000), no.~1, 1--15. 

 
\bibitem{TitsLNM}
J. Tits, {\it Buildings of spherical type and finite BN-pairs}, Lecture Notes in Mathematics, Vol. 386, Springer, Berlin, 1974. 

 
\bibitem{TitsFree}
J. Tits, Endliche Spiegelungsgruppen, die als Weylgruppen auftreten, Invent. Math. {\bf 43} (1977), no.~3, 283--295. 

 
\bibitem{TitsComo}
J. Tits, Immeubles de type affine, in {\it Buildings and the geometry of diagrams (Como, 1984)}, 159--190, Lecture Notes in Math., 1181 Springer, Berlin. 

 
\bibitem{WeSph}
R. M. Weiss, {\it The structure of spherical buildings}, Princeton Univ. Press, Princeton, NJ, 2003. 

\bibitem{WeAff}
R. M. Weiss, {\it The structure of affine buildings}, Annals of Mathematics Studies, 168, Princeton Univ. Press, Princeton, NJ, 2009. 



\end{thebibliography}
\end{document}